\documentclass[12pt]{article}

\usepackage{amsmath,amssymb,amsthm}
\usepackage{diagrams}
\chardef\bslash=`\\
\newtheorem{Thm}{Theorem}[section]
\newtheorem{Cor}[Thm]{Corollary}
\newtheorem{Lem}[Thm]{Lemma}
\newtheorem{Prop}[Thm]{Proposition}
\newtheorem{Def}{Definition}
\newtheorem{Rque}[Thm]{Remark}

\newtheorem{Exam}{Example}

\title{Poisson resolutions}
\author{Baohua Fu}

\def\cit{{\mathbb C}}

\def\qit{{\mathbb Q}}
\def\pit{{\mathbb P}}

\def\0{{\mathcal O}}

\def\T{{\mathcal T}}
\def\I{{\mathcal I}}
\begin{document}
\maketitle
\begin{abstract}
A resolution $Z \to X$ of a Poisson variety $X$ is called {\em Poisson} if every Poisson structure on $X$ lifts to a Poisson structure on $Z$.
For symplectic varieties, we prove that Poisson resolutions coincide with symplectic resolutions. It is shown that
for a Poisson surface $S$, the natural resolution $S^{[n]} \to S^{(n)}$ is a Poisson resolution. Furthermore, if
$Bs|-K_S| = \emptyset$, we prove that this is the unique projective Poisson resolution for $S^{(n)}$. 
\end{abstract}
\section{Definitions and introduction}

We work over the complex number field $\cit$. Let $X$ be an algebraic variety. A {\em Poisson structure} on $X$ is a skew-linear bracket
$\{, \}$ (called {\em Poisson bracket})
 on the structure sheaf $\0_X$ such that: (i) it is a derivation in each variable; (ii) it satisfies the Jacobi identity.
Every variety admits a trivial Poisson structure given by $\{ , \} \equiv 0.  $
A {\em Poisson variety} is an algebraic  variety which admits a non-trivial Poisson structure.

A Poisson bracket $\{ , \}$  on $X$ defines an $\0_X$-linear map  $\Theta:  \Omega_X^2 \to \0_X$ by $ \{f, g\} = \Theta(d f \wedge d g)$.
If we denote by $H_f$ the Hamiltonian vector field $\{ f, \cdot \}$, then the Jacobi identity reads $H_f(\Theta) = 0,$ i.e. Hamiltonian vector fields
preserve the Poisson structure.  

When $X$ is normal, giving a Poisson structure on $X$ is equivalent to giving a {\em Poisson bivector} 
 $\theta \in H^0(X_0, \wedge^2 \T_{X_0})$, where $X_0$ is the smooth part of $X$.
Then the Jacobi identity is equivalent to  $[\theta, \theta]=0,$ where $[  , ]$ is the usual Schouten bracket.
Note that giving  $\theta \in H^0(X_0, \wedge^2 \T_{X_0})$ is equivalent to giving a homomorphism of vector bundles
$B: T^*X_0 \to TX_0$, with $<B(\alpha), \beta> =  \theta (\alpha \wedge \beta). $ Then $[\theta, \theta]=0$ 
is equivalent to the vanishing of some quantity involving only $B$ (see Prop. 1.1 \cite{Bo1}).

 Any product of Poisson varieties is again Poisson.   The quotient of a Poisson variety $X$  by a finite group which preserves a Poisson structure
on $X$ is again a Poisson variety. It is proved by D. Kaledin that every irreducible component, 
any completion and the normalization of a Poisson variety is again Poisson (see \cite{Ka1}).

Any normal variety whose smooth part admits a  symplectic structure  is a Poisson variety.  This is the starting point for
the study of symplectic varieties by using techniques in   Poisson geometry, which is recently carried out by D. Kaledin (see \cite{Ka2}). 
This note can be regarded as a continuation in this direction. However, the point we want to make here is that the study of Poisson varieties 
and Poisson resolutions themselves may be of independent interest. Let us first 
recall some basic notions on symplectic varieties and symplectic resolutions.

Following A. Beauville (\cite{Bea}), a normal variety $W$ whose smooth part admits a 
symplectic structure $\omega$ is called a {\em symplectic variety} if 
$\omega$ lifts to a global 2-form $\tilde{\omega}$ (possibly degenerate) on any resolution of $W$. 
A resolution $\pi: Z \to W$ for a symplectic variety is called {\em symplectic} if 
$\tilde{\omega}$ is non-degenerate, i.e. if $\tilde{\omega}$ gives a symplectic structure on $Z$. 
It can be shown that $\pi$ is symplectic if and only if it is crepant.
This implies that $\pi$ is  a symplectic resolution if and only if 
for any symplectic structure $\omega'$ on $W_{reg}$, the lifted 2-form $\pi^* \omega'$ extends to 
a symplectic structure on the whole of $Z$. This motivates the following definition.

\begin{Def}
(i) Let $(X, \theta)$ be a Poisson variety endowed with a Poisson structure $\theta$. 
A {\bf Poisson resolution} for $(X, \theta)$ is a resolution  $\pi: Z \to X$ 
such that the Poisson structure $\theta$ lifts to a Poisson structure $\tilde{\theta}$ on $Z$.

(ii) A {\bf Poisson resolution} for a Poisson variety $X$ is a resolution
$\pi: Z \to X$ such that every Poisson structure on $X$ lifts to a Poisson structure on
$Z$ via $\pi$.

(iii) Two Poisson resolutions $\pi_i: Z_i \to X, i=1, 2,$ are called {\bf isomorphic} if
the rational map $\pi_2^{-1} \circ \pi_1: Z_1 --\to Z_2$ is an isomorphism.   
\end{Def}

It turns out that Poisson resolutions behave much like symplectic resolutions. The two notions coincide for 
symplectic varieties (Theorem \ref{sym}). Like that for symplectic resolutions, a resolution $\pi: Z \to X$ of
a Poisson variety $X$ with $codim(Exc(\pi)) \geq 2$ is always a Poisson resolution.
The uniqueness up to isomorphisms of symplectic resolutions for $S^{(n)}$ with $S$ symplectic can be generalized
to Poisson resolutions for $S$ Poisson with $Bs|-K_S| = \emptyset$ (Theorem \ref{unique}).
It would be interesting to exploit more such analogue properties, for example the semi-smallness.

However,  unlike that for symplectic resolutions, it is not clear that a product of Poisson resolutions is still
Poisson.  It would be interesting to find more relationships between crepant resolutions and Poisson resolutions.

Here is an outline of the contents:

$\bullet$ Section 2 is to study the liftability of Poisson structures to blowups;

$\bullet$ Section 3 studies Poisson resolutions for normal surfaces; 

$\bullet$ In Section 4, we prove that a resolution for a symplectic variety is Poisson if and only if it is symplectic;

$\bullet$ Section 5 proves that the natural resolution $S^{[n]} \to S^{(n)}$ is Poisson, where $S$ is  a smooth Poisson surface.
Furthermore, for surfaces $S$ such that $|-K_S|$ is base-point-free, we prove that this is the unique (up to isomorphisms) projective 
 Poisson resolution for $S^{(n)}$. \vspace{0.4cm}

{\em Acknowledgments:} The author wants to thank A. Beauville, M. Brion and P. Vanhaecke for helpful discussions.

\section{Blowups}
\begin{Prop}
Let $X$ be a smooth Poisson variety endowed with a Poisson structure $B: T^*X \to TX$ and $Y \subset X$ an irreducible smooth subvariety.
Consider the blowup of $X$ at $Y$: $ Z \xrightarrow{\pi} X$. If the Poisson structure $B$ lifts to $Z$, 
then $$N_{Y|X}^* \subseteq Ker B.$$ Furthermore, we have:

(i). when $codim Y = 2$, this condition is also sufficient;

(ii). when $Y$ is a point (denoted by $y$) and $dim X \geq 3$, then $B$ lifts to a Poisson structure on $Z$ if and only if:
 $B_y =0$ and $d (B(df)(g)) (y)= d (\theta (df \wedge dg)) (y) =0 $ for any local functions $f, g$ near $y$, where $\theta$
is the corresponding Poisson bivector on $X$. 
\end{Prop}
\begin{proof}
Take an open analytic subset $U \subset X$ and local coordinates $x_1, \cdots, x_n$ such that $Y \cap U = \{ x_{k+1} = \cdots = x_n =0 \}.$
Then $\pi^{-1}(U)$ can be covered by $n-k$ open sets, each having the following local coordinates: 
$$z_i = x_i,\ \text{for}\ i \in A \cup \{j\}; \quad  z_l = x_l / x_j, \ \text{for}\  l \in C,$$ where $k+1 \leq j \leq n$, 
  $A = \{1, \cdots, k\}$ and $C = \{k+1, \cdots, n\}- \{ j \}$.
Then $$ \partial_{x_i}  = \partial_{z_i},  \text{for}\ i \in A; \quad  \partial_{x_l}  = \frac{1}{z_j} \ \partial_{z_l},  \text{for}\ l \in C; \quad  \partial_{x_j} = \partial_{z_j} - \sum_{l \in C} \frac{z_l}{z_j}\ \partial_{z_l}.   $$

Write the Poisson bivector $\theta = \sum_{s < t} g_{st}(x) \partial_{x_s} \wedge  \partial_{x_t} .$
Set $$\tilde{g}_{lm} (z_1, \cdots, z_n) =g_{lm} (z_1, \cdots, z_k, z_jz_{k+1}, \cdots, z_j, \cdots, z_j z_n).$$
If $\theta$ lifts to a Poisson structure $\tilde{\theta}$ on $Z$, then 
$$\tilde{\theta} = H_1 + \frac{1}{z_j^2} H_2 + \sum_{i \in A, l\in C} \frac{\tilde{g}_{il} - z_l \tilde{g}_{ij}}{z_j} \partial_{z_i} \wedge \partial_{z_l} +\sum_{l\in C} \frac{\tilde{g}_{jl}}{z_j} \partial_{z_j} \wedge \partial_{z_l},  $$
where $H_1$ is holomorphic  and $H_2 = \sum_{l < m, l, m \in C} (\tilde{g}_{lm} - z_l \tilde{g}_{jm} + z_m \tilde{g}_{jl}) \partial_{z_l} \wedge  \partial_{z_m}$.

That $\tilde{\theta}$ is holomorphic implies that $\cfrac{\tilde{g}_{il} - z_l \tilde{g}_{ij}}{z_j}$ and $\cfrac{\tilde{g}_{jl}}{z_j} $ 
are both holomorphic. The latter is equivalent to $g_{il}(y) = g_{ij}(y) = g_{jl}(y) = 0 $ for any $y \in Y$, for any $i \in A, l\in C$ and for any $j$.
Notice that $B(d x_j) = \sum_i g_{i j} \partial_{z_{i}}, $
thus this is equivalent to that $B(d x_j) (y) = 0$ for $j = k+1, \cdots, n, y \in Y$, which shows that $N_{Y|X}^* \subseteq Ker B.$

Now if $codim\ Y = 2$, then $H_2 = 0$, thus this condition is also sufficient.

If $codim\ Y \geq 3$, then $\frac{1}{z_j^2} H_2$ is holomorphic if and only if $\frac{\partial H_2}{\partial z_j} $ is zero when $z_j =0 $ for any $j$.
The latter is equivalent to $\frac{\partial g_{lm}} {\partial x_s} |_Y = 0 $ for any $l, m, s \in C \cup \{j\},$ which can be reformulated as
 $[B(dx_l), \partial_{x_m}]|_Y = 0 $ for any $l, m \in C \cup \{j\}$. Now if $Y$ is just a point $y$, then this is equivalent to
that for any local functions $f, g, h$ near $y$,
$$0 = [B(df), \partial_h](g) (y) = B(df) (\partial_h g) (y) - \partial_h (B(df)(g)) (y) =  - \partial_h (B(df)(g)) (y),$$
where the last equality follows from $B(df) (y) = 0$. Notice that $\partial_h (B(df)(g)) (y) = 0  $ for any $h$
is equivalent to $d(B(df)(g)) (y) =0$, which completes the proof.

\end{proof}
\begin{Cor}\label{Cor22}
Let $S$ be a smooth surface equipped with a Poisson structure $\theta \in H^0(S, \omega_S^{-1})$. Take a point $x \in S$ and
let $S_x \to S$ be the blowup of $S$ at $x$.  Then the Poisson structure $\theta$ lifts to $S_x$ if and only if $\theta(x) =0$.  
\end{Cor}
\begin{Rque}
The condition that $X$ is smooth can be replaced by the condition that $Y$ is contained in the smooth part of $X$.
\end{Rque}


\begin{Exam}\label{exam1}
Let $S$ be a smooth surface endowed with a Poisson structure $s \in H^0(S, \omega_S^{-1}).$ Suppose that there exists a smooth
curve $C$ contained in the support of $div(s)$. Take an element $\I \in S^{[2]}$, then $T^*_{\I}(S^{[2]})$ (resp. $T_{\I}(S^{[2]})$)
is isomorphic to $Hom(\I, \0_S/\I \otimes \omega_S)$ (resp.  $Hom(\I, \0_S/\I)$). Now by the multiplication by $s$, we obtain
a linear map $T^*_{\I}(S^{[2]}) \to T_\I(S^{[2]})$, which gives a Poisson structure $\theta$ on $S^{[2]}$ (for details, see \cite{Bo2}).

As easily seen, the Poisson structure is zero when restricted to $C^{(2)} \subset S^{[2]}$. If we denote by $Z$ the blowup of
$S^{[2]}$ at $C^{(2)}$, then by the precedent proposition, the Poisson structure $\theta$ lifts to a Poisson structure on $Z$.
\end{Exam}

\section{Normal surfaces}   
Let $S$ be a normal surface with $S_0$ its smooth part. The closure in $S$ of a canonical divisor of $S_0$ defines a Weil divisor (class)
$K_S$ on $S$, which is called the canonical divisor of $S$.  The dualizing sheaf $\omega_S$ is the sheaf $\0_S(K_S)$, which is nothing
but $j_* \omega_{S_0}$, where $j: S_0 \to S$ is the natural inclusion. We denote by $\omega_S^{-1}$ the sheaf $\0_S(- K_S)$.
\begin{Prop}
$S$ carries a (non-trivial) Poisson structure if and only if $H^0(S, \omega_S^{-1}) \neq 0$. 
\end{Prop}
\begin{proof}
Since $S$ is normal, $Sing(S)$ consists of points, thus  $S$ admits   a Poisson structure if and only if $S_0$ admits a 
Poisson structure. The latter is equivalent to $H^0(S_0, \omega_{S_0}^{-1}) \neq 0$. By definition,
$H^0(S, \omega_S^{-1}) = H^0(S,j_* \omega_{S_0}^{-1}) = H^0(S_0, \omega_{S_0}^{-1}), $ which completes the proof.
\end{proof}

Let $D$ be a Weil divisor on $S$ and $\pi: \tilde{S} \to S$ a resolution. We denote by $E_i, i=1, \cdots, k$ the irreducible
components of $Exc(\pi)$. The pull-back of $D$ is defined by $\pi^*(D) = \bar{D} + \sum_i a_i E_i$, where
$\bar{D}$ is the strict transform of $D$ and $a_i, i=1, \cdots, k$ are rational numbers determined by $\pi^*(D) \cdot E_i = 0$.
Furthermore if $D$ is effective, then  $\pi^*D$ is effective (see \cite{Mum}).

Now let $Z = \pi^* K_S - K_{\tilde{S}}$, then $Z$ is a $\qit$-divisor with support in  $Exc(\pi)$. 
$Z$ is called the {\em canonical cycle} of the resolution $\pi$. If we write $Z = \sum_i b_i E_i$, then the numbers
$b_i, i=1, \cdots, k$ are determined by $Z \cdot E_i = - K_{\tilde{S}} \cdot E_i. $

\begin{Lem} \label{iff}
A resolution $\pi: \tilde{S} \to S$  is   Poisson  if and only if for any $F \in |-K_S|$, the divisor $ \pi^* F + Z$
is effective, where $Z$ is the canonical cycle of $\pi$.
\end{Lem}
\begin{proof}
Let $\theta \in H^0(S, \omega_{S}^{-1})$ be the section (up to scalars) defined by $F$. If $\theta$ lifts to  a section $\tilde{\theta}$ 
of the sheaf $\omega_{\tilde{S}}^{-1}$,  then $div(\tilde{\theta}) = \pi^* F + Z$ is effective. Conversely, if  $\pi^* F + Z$
is effective, then it gives the section $\tilde{\theta}$ lifting $\theta$. Thus $\pi$ is a Poisson resolution if and only if
for any $F \in |-K_S|$, the divisor $\pi^* F + Z$ is effective.
\end{proof}

\begin{Prop}
 A minimal resolution is a Poisson resolution.
\end{Prop}
This follows from the fact that the canonical cycle of a minimal resolution is effective. Any other Poisson resolution for a normal
surface can be otained from some blowups of the minimal resolution. The condition for the blowup center is given by 
 Corollary \ref{Cor22}.




\begin{Prop}
Let $S$ be a Poisson normal surface. Then $S$ admits a unique (up to isomorphisms) Poisson resolution (given by the minimal resolution)
if and only if: (i) $S$ has at worst RDPs as singularities; (ii)   $Bs|-K_S| \cap Sing(S) = \emptyset$.
\end{Prop}
\begin{proof}
Consider the minimal resolution $\pi: S_{min} \to S$. Notice that the canonical cycle $Z$ is contained in $Bs |-K_{S_{min}}|$.
By Corollary \ref{Cor22}, if $S$ admits a unique Poisson resolution, then $Z =0$, which shows (i). Now suppose that
$Bs|-K_S| \cap Sing(S) \neq \emptyset$. Take a point $x \in Bs|-K_S| \cap Sing(S)$, then for any $F \in |-K_S|$, the 
vector $v = (-\bar{F} \cdot E_1, \cdots, - \bar{F} \cdot E_k) $ is a non-zero vector with each coordinate $\leq 0$,
where $E_i$ are irreducible components of $\pi^{-1}(x)$. Then $\pi^*(F) = \bar{F} + \sum_i a_i E_i$, where
$(a_1, \cdots, a_k) = v M^{-1}$ and $M = (E_i \cdot E_j)_{k \times k}.$  By (i), $x$ is a RDP, thus $M$ is a Cartan matrix of type
ADE, thus $M^{-1}$ is a matrix with strictly negative entries, which shows that $a_i > 0$ for any $i$.
This proves that $\pi^{-1}(x) \subset Bs|-K_{S_{min}}|$. Now by Corollary \ref{Cor22}, any blowup
of $S_{min}$ at a point in  $\pi^{-1}(x)$ is again a Poisson resolution for $S$, thus   $S$ admits many Poisson resolutions.

Now  suppose that $S$ has only RDPs and  $Bs|-K_S| \cap Sing(S) = \emptyset$. 
Then for any $x \in Sing(S)$, there exists an element
$F \in Bs|-K_S|$ such that $x \notin F$. Then $supp(\pi^* F) \cap \pi^{-1}(x) = \emptyset.$ This implies that
$Bs |-K_{S_{min}}| \cap \pi^{-1}(Sing (S)) = \emptyset.$ Now by Corollary \ref{Cor22}, $S$ admits a unique Poisson resolution.
\end{proof}

Recall that a Gorenstein log Del Pezzo surface is a  projective normal surface with only RDPs and with an ample anti-canonical divisor.
\begin{Cor}
Let $S$ be a symplectic surface or a Gorenstein log Del Pezzo surface.  Then $S$ admits a unique Poisson resolution.
\end{Cor}

\begin{Exam} \label{exam}
 Consider the quartic surface $S =\{ x_0x_1^3 + x_0x_2^3+x_3^4= 0\} \subset \pit^3$.
$S$ has trivial dualizing sheaf, thus its smooth part $S_0$ admits a symplectic structure $\omega$, but $S$ is not a symplectic surface,
 since $Sing(S) = \{[1:0:0:0]\}$  is a minimally elliptic singularity. Thus
 for any resolution $\pi: \tilde{S} \to S$, $\omega$ does not lift to a 
global 2-form on $\tilde{S}$. However, when $\pi$ is minimal, $\omega$ lifts to a Poisson structure on $\tilde{S}$.
\end{Exam}

\section{Symplectic varieties}
\begin{Thm}\label{sym}
Let $W$ be a symplectic variety. Then a resolution $\pi: Z \to W$ is a  symplectic resolution if and only if it is a Poisson resolution.
\end{Thm}
\begin{proof}
Let $\pi: Z \to W$ be a symplectic resolution with $\tilde{\omega}$ (resp. $\omega$) the symplectic form on
$Z$ (resp. $W_0$), where $W_0$ is the smooth part of $W$.
 Giving  a Poisson structure $\theta$ on $W$ is equivalent to giving a Poisson structure on
$W_0$. Now the symplectic structure $\omega$ gives an isomorphism between $\Omega^2_{W_0}$ and $\wedge^2T_{W_0}$,
thus  $\theta$ corresponds to a 2-form $\alpha$ on $W_0$. We claim that $\pi^* \alpha$ extends to a 2-form $\tilde{\alpha}$
on the whole of $Z$.

Notice that the extension is unique if exists, thus we need only to show this claim locally.
In fact, for any $x \in W$, we take an affine open set $x \in U \subset W$, then $U$ is Stein.  Now by Theorem 4. ([Nam])
we have $\pi_* \Omega^2_{\pi^{-1}(U)} = i_* \Omega^2_{U_0}$, where $i: U \to W$ is the inclusion and $U_0$ is the smooth part of $U$. This shows that
$\pi^*(\alpha|_{U_0})$ can be extended to $\pi^{-1}(U)$. 

In conclusion, we get a global 2-form $\tilde{\alpha}$ on $Z$, which gives a bivector $\tilde{\theta}$ on $Z$ via the symplectic
form $\tilde{\omega}$.  Now $[\tilde{\theta}, \tilde{\theta}] =0$ follows from $[\theta, \theta]=0$, thus
$\tilde{\theta}$ gives a Poisson structure on $Z$ lifting the Poisson structure $\theta$, which shows that $\pi$ is a Poisson resolution.

Now suppose that $\pi$ is a Poisson resolution. Consider the Poisson structure on $W_0$ defined by the symplectic structure
$\omega$. Let $Z_0 =\pi^{-1}(W_0)$. Then $\omega$ defines a homomorphism of vector bundles $B_0: T_{Z_0}^* \to T_{Z_0}$
by $\pi^*(\omega) (u \wedge B_0(\alpha)) = \alpha(u),$ which has an inverse $C_0:  T_{Z_0} \to T_{Z_0}^*$ defined by
$C_0(u)(v) = \pi^*(\omega)(v \wedge u).$  By the definition of symplectic varieties, $\pi^*\omega$ extends to a 2-form $\tilde{\omega}$
on the whole of $Z$. This gives an extension of $C_0$ to a homomorphism $C: T_Z \to T_Z^*$ by $C(u)(v) = \tilde{\omega}(v \wedge u).$
Notice that $C$ is invertible if and only if $\tilde{\omega}$ is non-degenerate.

Now $\pi$ is Poisson implies that $B_0$ extends to a homomorphism $B: T_{Z}^* \to T_{Z}$. Furthermore this homomorphism
satisfies $\tilde{\omega} (u \wedge B(\alpha)) = \alpha(u).  $ Now consider the composition $CB: T^*_Z \to T_Z^*$, which verifies
$$CB(\alpha)(v) = \tilde{\omega} (v \wedge B(\alpha)) = \alpha(v),    $$
for any local vector field $v$. 
Thus $CB = id.$ In particular, $C$ and $B$ are bother invertible, which shows that  $\tilde{\omega}$ is symplectic on the whole
of $Z$, thus $\pi$ is a symplectic resolution.
\end{proof}
\begin{Rque}
In the proof above, we used in an essential way that $\pi^*\omega$ extends to a global 2-form on $Z$. In fact, let $S$ be the quartic
surface in Example \ref{exam}. Then $S$ admits a symplectic structure, but it is not a symplectic variety. $S$ admits many Poisson resolutions, but
non of them is symplectic. However, we have the following easy corollary.
\end{Rque}
\begin{Cor}
Let $W$ be a normal variety whose smooth part admits a symplectic structure. Then any crepant resolution for $W$ is a Poisson resolution.
\end{Cor}
\begin{proof}
Let $\pi: Z \to W$ be  a crepant resolution for $W$.  Since the smooth part of $W$ carries a symplectic structure, then $K_W$ and $K_Z$ are
both trivial. Now a classical result of Elkik and Flenner implies that $W$ has only rational singularities. Now by \cite{Nam}, 
$W$ is a symplectic variety, so $\pi$ is a symplectic resolution, which is thus a Poisson resolution. 
\end{proof}

It would be interesting to construct an example (if there exists one) of a crepant resolution for a Poisson variety which is not a Poisson resolution. 

For other interactions  between Poisson structures and symplectic varieties, the interested reader is refereed to D. Kaledin's paper \cite{Ka2}. 
\section{Hilbert schemes}
\begin{Thm} 
Let $S$ be a smooth Poisson surface and $\pi: S^{[n]} \to S^{(n)}$ the natural resolution of the $n$th symmetric product $S^{(n)}$.
Let $W$ be a symplectic variety which admits a symplectic resolution: $Z \to W$. Then

(i). the resolution $\pi: Z \times S^{[n]} \to W \times S^{(n)}$ is a Poisson resolution;

(ii). any crepant resolution for $W \times S^{(n)}$ is a Poisson resolution.
\end{Thm}
\begin{proof}
Let $X = Spec(\cit[x,y,z]/(xy = z^2))$ and $ \tilde{X} \to X$ the blowup at the origin. 
Then $Z \times  \cit^{2n-2} \times \tilde{X} \to  W \times \cit^{2n-2} \times X$ is a symplectic resolution, thus it is a Poisson resolution,
by Theorem \ref{sym}.  

Now let $\Delta = \{2x_1 + x_2 + \cdots +x_{n-1} | x_i \neq x_j\} \subset S^{(n)}$ and $S_*^{(n)} = \Delta \cup (S^{(n)} - Sing (S^{(n)})). $
Set $S_*^{[n]} = \pi^{-1} (S_*^{(n)}).$ Then locally  $Z \times S_*^{[n]} \to W \times S_*^{(n)}  $ 
 is isomorphic to the resolution 
$Z \times \cit^{2n-2} \times \tilde{X} \to W \times \cit^{2n-2} \times X$, 
thus every Poisson structure $\theta$ on $W \times S^{(n)}$ can be lifted to a Poisson structure $\tilde{\theta}$
on  $Z \times S_*^{[n]}$. Notice that the 
complement of $Z \times S_*^{[n]}$ in $Z \times S^{[n]}$ has codimension 2, thus $\tilde{\theta}$ extends to the whole of $Z \times S^{[n]}$.  

Now suppose that we have another crepant resolution $p: Y \to W \times S^{(n)}$, then the rational map
$\phi= \pi^{-1} \circ p: Y --\to Z \times S^{[n]}$ is isomorphic in codimension 1, since $\pi$ and $p$ are both crepant.
So every Poisson structure on $Z \times S^{[n]}$ induces a Poisson structure on $Y$ through $\phi$, thus $p$ is
also a Poisson resolution. 
\end{proof}
\begin{Cor} \label{hilb}
Let $S$ be a smooth Poisson surface. Then the natural resolution $S^{[n]} \to S^{(n)}$ is a Poisson resolution.
\end{Cor}
\begin{Rque}
Every Poisson structure on $S$ defines a Poisson structure on $S^{(n)}$.  
It is proved by F. Bottacin (\cite{Bo2}) that these Poisson structures lift to $S^{[n]}$. He 
explicitly constructed the lifted Poisson structures (see Example \ref{exam1}). 
\end{Rque}
\begin{Rque}
The assumption $S$ being   Poisson is not necessary, what we need is  that $S^{(n)}$ is a Poisson variety.  Another
proof of the theorem goes as follows: Let $E$ be the exceptional divisor of $S^{[n]} \to S^{(n)}$. 
Take a Poisson structure $s \in H^0(S, \omega_S^{-1})$
and let $F = div(s)$. Then $s$ induces a symplectic structure on $U = S - supp(F)$.  Then $Z \times U^{[n]} \to W \times  U^{(n)}$ 
is a symplectic resolution, thus every Poisson structure $\theta$ on $W \times S^{(n)}$ extends to $Z \times U^{[n]}$.
Notice that $(Z \times E) \cap (Z \times U^{[n]})$ is open, thus $\theta$ extends to the generic point of  $Z \times E$. This shows that
$\theta$ extends to the whole of $Z \times S^{[n]}$.  This argument can be applied to prove the following: 
{\em Let $S_i$ be smooth Poisson surfaces and $Z \to W$ a symplectic resolution, then $Z \times \prod_i S_i^{[n_i]} \to W \times \prod_i S_i^{(n_i)}$ is a Poisson resolution. }
\end{Rque}

We should point out  that 
$W \times S^{(n)}$ can have more than one crepant resolutions, since $W$ may have many symplectic resolutions.
However,  we do not know if every crepant resolution for $W \times S^{(n)}$
is a product of $S^{[n]}$ and a symplectic variety. This is true if $W$ is itself smooth, by Theorem 2.2 \cite{FN} (see
also Remark \ref{rque}). It looks like that every Poisson resolution for $W \times S^{(n)}$ is crepant, however, we do not know how to
prove it.

 Now we come to study the uniqueness of Poisson resolutions in some special cases.

\begin{Thm} \label{unique}
Let $S$ be a smooth Poisson surface such that $|-K_S|$ is base-point-free.
Then any projective Poisson resolution for $S^{(n)}$ is isomorphic to the natural resolution: $S^{[n]} \xrightarrow{\pi} S^{(n)}$.
\end{Thm}
\begin{proof}
Take a  projective Poisson resolution $\pi_1: Z \to  S^{(n)}$ and let $\phi = \pi^{-1} \circ \pi_1$ be the rational map $Z--\to S^{[n]}$ .
We need to prove that $\phi$ is an isomorphism.  Since $ S^{(n)}$ is normal and $\qit$-factorial,  the exceptional locus 
$E_1 = exc(\pi_1)$ is of pure codimension 1. It is well-known that the exceptional locus of $\pi$ (denoted by $E$) is an irreducible
divisor. The main difficulty is to prove that $E_1$ is also irreducible. 

Take a point $I = \sum_i k_i x_i$ on $S^{(n)}$. 
By the hypothesis that $Bs|-K_S| = \emptyset$, there exists an element $F \in |-K_S|$ such that $U: = S - supp(F)$
contains all $x_i$.  Then $I \in U^{(n)}$. Notice that
$F$ defines a Poisson structure $s$ (up to scalars)  on $S$, which induces a symplectic structure $\omega$  on $U$. Now $s$ defines
a Poisson structure $\theta$ on $S^{(n)}$, which when restricted to $U^{(n)}$ is the Poisson structure coming from $\omega$. 
 Note that $\pi_1$ lifts $\theta$
to a Poisson structure on $Z$, by Theorem \ref{sym}, $\pi_1^{-1} (U^{(n)}) \xrightarrow{\pi_1|_U} U^{(n)}$ is a projective symplectic resolution.
By Theorem 2.2 (\cite{FN}), this is isomorphic to the natural resolution $U^{[n]} \to U^{(n)}.$ In particular, the exceptional locus
$E_U: = exc(\pi_1|_U)$ is an irreducible divisor in $\pi_1^{-1}(U^{(n)})$.

Now for any other point  $I' = \sum_i k'_i y_i$ on $S^{(n)}$. The same argument gives another open set $V \subset S$ such that
$\pi_1^{-1}(V^{(n)}) \xrightarrow{\pi_1|_V} V^{(n)}$ is isomorphic to $V^{[n]} \to V^{(n)},$ and $I' \in V^{(n)}$. Thus $E_V: = exc(\pi_1|_V)$
is again an irreducible divisor. Notice that $E_V \cap E_U$ is open, thus $\overline{E_V} = \overline{E_U}$ on $Z$. Since this is true for
any point on $S^{(n)}$, the divisor $E_1$ is irreducible. Furthermore, this shows that $\phi: Z --\to S^{[n]}$ is isomorphic in codimension 1.

Now a general argument shows that $-E_1$ is $\pi_1$-ample, whose strict transform by $\phi$ is $-E$. The latter being $\pi$-ample implies
that $\phi$ is in fact an isomorphism. 
\end{proof}

\begin{Rque}
For smooth Poisson surfaces $S$  such that $Bs|-K_S| = \emptyset$, the precedent proof shows that the resolution $\pi: S^{[n]} \to S^{(n)}$
is locally symplectic resolutions, thus by Theorem \ref{sym}, $\pi$ is a Poisson resolution. This gives another proof of Corollary \ref{hilb}
in this special case.
\end{Rque}

\begin{Rque}
The condition that $S$ is smooth cannot be weakened. In fact if we take $S = Spec(\cit[x,y,z]/(xy = z^2))$ and $p: S_m \to S$ the minimal resolution.
Put $C = exc(p)$, which is $\pit^1$. Now $S^{(2)}$ admits two non-isomorphic Poisson resolutions (which are also symplectic)
given by $S_m^{[2]}$ and the Mukai flop along $C^{(2)} \subset S_m^{[2]}.$
\end{Rque}

\begin{Rque} \label{rque}
The above argument shows in fact the following stronger result: Let $S$ be a smooth Poisson surface with $Bs|-K_S| = \emptyset$. Let $X$
be a smooth symplectic variety. Then any projective Poisson resolution for $X \times S^{(n)}$ is isomorphic to
 $X \times S^{[n]} \to  X \times S^{(n)}$.
\end{Rque}

Using a similar argument as above, and that in the proof of Theorem 2.2 \cite{FN}, we can prove the following, which
 gives a generalization of Corollary \ref{hilb} and Theorem \ref{unique}.
\begin{Thm}
Let $S_i, i=1, \cdots, k$  be smooth Poisson surfaces  such that $|-K_{S_i}|$ is base-point-free. Then the Poisson variety $\prod_{i=1}^k S_i^{(n_i)}  $
admits a unique projective Poisson resolution (up to isomorphisms) given by: $$ \prod_{i=1}^k S_i^{[n_i]} \to  \prod_{i=1}^k S_i^{(n_i)}.  $$
\end{Thm}

Now we end this section by a classification of projective smooth Poisson surfaces $S$ with $Bs|-K_S| = \emptyset$,
which could be of independent interest.  The proof is 
a case-by-case check based on the classification (see \cite{BM}) of projective smooth Poisson surfaces. 
Notice that there are many other quasi-projective
smooth Poisson surfaces with $Bs|-K_S| = \emptyset$.
\begin{Prop}
Let $S$ be a projective smooth Poisson surface such that $|-K_S|$ is base-point-free. Then $S$ is isomorphic to one of the following surfaces:

(i).  $K3$ surfaces or abelian surfaces;

(ii). blowup of $k \leq 8$ points in general position in $\pit^2$; 

(iii). blowup of $k \leq 7$ points in general position in $\pit^1 \times \pit^1$ or in $\pit(\0_{\pit^1} \oplus \0_{\pit^1}(2))$;

(iv). blowup of $k \leq 1$ general point in $\pit^1 \times E$, where $E$ is an elliptic curve. 
\end{Prop}

\quad \\
Labortoire J. A. Dieudonn\'e, Parc Valrose \\ 06108 Nice cedex 02, FRANCE \\
baohua.fu@polytechnique.org
\end{document}